\documentclass[reqno,12pt]{amsart}
\usepackage{a4wide}
\usepackage{amscd}
\input{epsf}

\numberwithin{equation}{section}

\newcommand\ce{\ensuremath{\mathfrak{ e}}}
\newcommand\cf{\ensuremath{\mathfrak{ f}}}

\newcommand\co{\ensuremath{\mathfrak{ o}}}
\newcommand\cp{\ensuremath{\mathfrak{ p}}}

\newcommand\cs{\ensuremath{\mathfrak{ s}}}

\newcommand\cu{\ensuremath{\mathfrak{ u}}}

\newcommand\cG{\ensuremath{\mathfrak{ G}}}

\newcommand\cL{\ensuremath{\mathfrak{ L}}}

\newcommand\cP{\ensuremath{\mathfrak{ P}}}

\newcommand\tr{\ensuremath{\mathrm{ tr\,}}}
\newcommand\grad{\ensuremath{\mathrm{ grad}}}
\newcommand\Aut{\ensuremath{\mathrm{ Aut}}}
\newcommand\Diag{\ensuremath{\mathrm{ Diag}}}
\renewcommand\Im{\ensuremath{\mathrm{ Im}}}
\renewcommand\Re{\ensuremath{\mathrm{ Re}}}

\begin{document}

\title{Projective Planes, Severi Varieties and Spheres}
\author{Michael Atiyah and J\"{u}rgen Berndt}
\address{University of Edinburgh\\ Department of Mathematics and
Statistics\\ Mayfield Road\\
Edinburgh EH9~3JZ\\ United Kingdom}
\email{atiyah@maths.ed.ac.uk}
\address{University of Hull\\ Department of Mathematics\\ Cottingham Road\\ Hull
HU6 7RX\\ United Kingdom}
\email{j.berndt@hull.ac.uk}
%\subjclass[2000]{...}
%\date{June 5, 2002}

\begin{abstract}
A classical result asserts that the complex projective plane modulo
complex conjugation is the $4$-dimensional sphere. We generalize
this result in two directions by considering the
projective planes over the normed real division algebras and by
considering the complexifications of these four projective planes.
\end{abstract}

\maketitle
\thispagestyle{empty}

\bigskip\bigskip\bigskip
\section{Introduction}

There is an elementary but very striking result which asserts that
the quotient of the complex projective plane $CP^2$ by complex
conjugation is the $4$-dimensional sphere. This result has
attracted the attention of many geometers over the years,
rediscovered afresh each time and with a variety of proofs. For
some historical remarks about the origins of this theorem we refer
to Arnold \cite{Arn}, where it is put in a more general context.

The purpose of the present paper is to give two parallel
generalizations of this theorem. In the first we view $CP^2$ as
the second member of the family of projective planes over the four
normed real division algebras. This is close to Arnold's treatment
and completes it by dealing with the octonions. The second
generalization views $CP^2$ as the complex algebraic variety
obtained by complexifying the real projective plane $RP^2$ and
hence as the first member of the four algebraic varieties got by
complexifying the four projective planes. This family of algebraic
varieties has appeared in several contexts. First, in relation to
Lie groups and the ``magic square'' of Freudenthal \cite{Fre1},
somewhat independently in the characterization by Zak \cite{Zak}
of ``Severi varieties'', and also in the classification by
Nakagawa and Takagi \cite{NaTa} of K\"ahler submanifolds with
parallel second fundamental form in complex projective spaces.
Since these varieties are not so widely known we shall give a
brief account of them in an appendix. An extensive treatment based
on the ``magic square'' can be found in \cite{LaMa}.

Denoting by $CP^2$, $HP^2$ and $OP^2$ the projective planes over
the complex numbers $C$, the quaternions $H$ and the octonions
$O$, and by $S^d$ the sphere of dimension $d$, our first result
may be formulated as follows.

\medskip
{\it There are natural diffeomorphisms}
\begin{equation}
\begin{array}{rcl}
CP^2/O(1) & = & S^4\ ,\\ HP^2/U(1) & = & S^7\ , \\ OP^2/Sp(1) & =
& S^{13}\ .
\end{array} \label{quotient}
\end{equation}

\medskip
{\sc Remarks.} 1) The first equality is the ``folklore'' theorem
which provides our starting point. The second one is proved in
\cite{Arn} and independently in \cite{AtWi}.

2) The maps from the projective planes to the spheres in
(\ref{quotient}) are fibrations outside the ``branch locus'' given
by the preceding projective plane. The sense in which the
equations
in (\ref{quotient}) are diffeomorphisms is explained in the next
section.

3) The embeddings of the branch loci
\begin{equation}
RP^2 \subset S^4\ ,\ CP^2 \subset S^7\ ,\ HP^2 \subset S^{13}
\label{branchloci}
\end{equation}
are well-known and will be elaborated on in Section 3.

4) We shall later in Section 4 formulate (\ref{quotient}) more
precisely as Theorem A. This will include an explicit and simple
construction of the maps from the projective planes to the
spheres. These maps will also be compatible with the relevant
symmetry group $SO(3)$, $SU(3)$ and $Sp(3)$.

\medskip
Before formulating our second result we need to introduce the {\it
complexifications} of the four projective planes. The first case
is clear, it leads from $RP^2$ to $CP^2$. Note that the action of
$SO(3)$ on $RP^2$ extends to a complex action of $SO(3,C)$ on
$CP^2$ and this leaves invariant the complex curve $z_1^2 + z_2^2
+ z_3^2 = 0$ which, for reasons that will be clear later, we
denote by $CP^2(\infty)$. All this generalizes to the other
projective planes. If we denote the projective planes by $P_n$ ($n
= 0,1,2,3$), so that $\dim P_n = 2^{n+1}$, then their
complexifications\footnote{$P_n(C)$ should not be confused with
the complex projective space of dimension $n$, for which we use
the notation $CP^n$.} $P_n(C)$ are complex algebraic varieties of
complex dimension $2^{n+1}$. The isometry group of $P_n$ extends
to an action of its complexification on $P_n(C)$ leaving invariant
a complex hypersurface $P_n(\infty)$. Moreover $P_n(C)$ has a
``real structure'', i.e.\ a complex conjugation. The real points
are just $P_n$, and $P_n(\infty)$ inherits a real structure with
no real points. We can now state our second result.

\medskip
{\it For each $n = 0,1,2,3$ we have a natural map
\begin{equation}
\varphi_n : P_n(C) \to S^{d(n)}\ ,\ d(n) = 3\cdot 2^n + 1\ ,
\label{mapphin}
\end{equation}
which is a fibration outside the branch locus $P_n$ and the
hypersurface $P_n(\infty)$. The fibres are the spheres
$S^0,S^1,S^3,S^7$.}

\medskip
{\sc Remarks}. 1) The case $n = 0$ of (\ref{mapphin}) is just the
first case of (\ref{quotient}).

2) The image of $P_n(\infty)$ under $\varphi_n$ is the focal set
of the branch locus $P_n$ in $S^{d(n)}$ and turns out to be
another embedding of $P_n$ in $S^{d(n)}$: it is just the image
under the antipodal map of the embedding $\varphi_n(P_n)$. For $n
< 3$ this is just the embedding in (\ref{branchloci}), and for $n
= 3$ it is the known fourth member of this sequence (as will be
explained in Section 4).

3) The varieties $P_n(C)$ are in fact homogeneous spaces of a
larger group as we shall explain later. For $n = 1$
the complexification of $CP^2$ is just the product $CP^2 \times
CP^2$, while for $n = 2$ the complexification of $HP^2$ turns out
to be the complex Grassmannian $Gr_2(C^6)$ of $2$-planes in $C^6$,
and for $n=3$ the complexification of $OP^2$
is the exceptional Hermitian symmetric space
$E_6/Spin(10)U(1)$.

\medskip
The paper is organized as follows. In Section 2 we review
well-known elementary properties of the normed real division
algebras and the associated geometries. In particular we explain
the general notion of a branched fibration of which
(\ref{mapphin}) is an illustration. In Section 3 we move on to
projective planes and their automorphisms, where we pay special
attention to the Cayley plane and its exceptional status. Section
4 examines the orbit structure for the projective planes and
spheres in (\ref{quotient}) with respect to the relevant symmetry
groups. This enables us to establish (\ref{quotient}). Section 5
is devoted to the formulation and proof of Theorem A, the more
precise version of (\ref{quotient}), while Section 6 deals
similiarly with Theorem B, a more explicit version of
(\ref{mapphin}). In Section 7 we establish projective versions of
the two theorems, where the symmetry group is replaced by the
larger group of projectivities of the relevant projective planes.
In the appendix we elaborate more on the
complexified projective planes $P_n(C)$.

\medskip
We thank Friedrich Hirzebruch and Jean-Pierre Serre for drawing
our attention to some of the relevant literature.

\section{Projective Lines, Hopf Maps and Branching}

We recall that there are four normed division algebras $A_n$ over
the reals, with $\dim A_n = 2^n$, namely

\smallskip
$n=0$\ : $R$, the real field, of dimension $1$;

$n=1$\ : $C$, the complex field, of dimension $2$;

$n=2$\ : $H$, the non-commutative field of quaternions, of dimension
$4$;

$n=3$\ : $O$, the non-associative algebra of octonions, of dimension
$8$.

\smallskip
To each of these we can associate the corresponding projective
line $A_nP^1$ by adding a point $\infty$ to the algebra. This
shows that the projective line is a sphere. The more classical
definition of a projective line is to consider all lines (i.e.\
one-dimensional subspaces) in the $2$-dimensional vector space
$A_n^2$ over $A_n$. For $n=3$ one has to be a little careful with the
definition of a line because of the non-associativity of the octonions.
A comprehensive introduction to the octonions can be found in
the survey article \cite{Bae}.

Now consider the tautological line bundle over
$A_nP^1$, where the fibre over each point in $A_nP^1$ is the
corresponding line in $A_n^2$.
It can also
be obtained by cutting the sphere $A_nP^1$ into two closed discs and
then gluing together the trivial line bundles over these two discs
along the boundary by using the multiplication with
elements of norm one in the
normed division algebra $A_n$.
We denote by $\cL_n$ the dual
bundle of the tautological line bundle.
The line bundle $\cL_n$ generates the reduced real K-theory of the
sphere of dimension $2^n$, and the octonionic line bundle $\cL_3$
induces a periodicity between the reduced real K-theories of
higher-dimensional spheres, known as Bott periodicity.

Since each fibre of $\cL_n$
is equipped with a
norm this line bundle naturally induces a sphere bundle $S(\cL_n)$ over
$A_nP^1$. The manifolds $S(\cL_n)$ are again spheres.
The fibrations $S(\cL_n) \to A_nP^1$ are
usually called the Hopf fibrations (though for $n=1$ it goes back originally
to W.K.~Clifford). In detail they are
\begin{equation}
\begin{array}{rcl}
\CD S^1 @>O(1)>>RP^1 &\ = &\  S^1\ ,\\
S^3 @>U(1)>>CP^1 &\ = & \ S^2\ ,\\
S^7 @>Sp(1)>>HP^1 &\ = & \ S^4\ ,\\
S^{15} @>S^7>>OP^1  &\ = &\ S^8\ . \endCD
\end{array}
\label{Hopffibration}
\end{equation}

\medskip
\noindent
The first three fibrations are principal (i.e.\ group actions)
while the last is not:
$S^7$ is the set of elements of norm one in the octonions but is not a
group since $O$ is non-associative.
The fact that every division algebra over
$R$ has dimension $2^n$, $n=0,1,2,3$,
can be proved topologically by showing that there are no
further sphere fibrations beyond (\ref{Hopffibration}).
More details about the construction of these fibrations
can be found in \S 20 of \cite{Ste}.
The fibres in
(\ref{Hopffibration}) will be denoted by $\Gamma_n$,
so that
\begin{equation}
\begin{array}{lllll}
\Gamma_0 & = & S^0 & = & O(1) \ , \\
\Gamma_1 & = & S^1 & = & U(1) \ , \\
\Gamma_2 & = & S^3 & = & Sp(1)\ , \\
\Gamma_3 & = & S^7 & . &
\end{array}
\label{fibres}
\end{equation}
We will discuss the Hopf fibrations again in Section 3 in relation
with group actions.

We now discuss the notion of ``branching'' as we shall encounter
it in this paper. The classical situation occurs in complex
variable theory where one Riemann surface can appear as the
branched covering of another. We shall restrict ourselves to the
simple case of double coverings, where the local model is the
equation $w = z^2$. Although the group of order $2$, given by $z
\mapsto -z$, has a fixed point at $z = 0$, the quotient is still a
smooth surface. The underlying topological reason is that, on the
small circles ${\mid} z {\mid} = \epsilon$
surrounding the fixed point,
we get the double covering in the first line of
(\ref{Hopffibration}) so that the quotient is still a circle and
hence is the boundary of a small disc.

If we take the product with $R^{n-2}$ we get the more general
situation where a group of order $2$ acting on an $n$-dimensional
manifold, with fixed-point components all of codimension $2$, has
a manifold as quotient. Again we refer to the fixed-point set as
the branch locus of the double covering.

The purpose of this lengthy analysis of a familiar situation was
to point out that each of the equations in (\ref{Hopffibration})
gives rise to a similar story, except that the finite group $O(1)$
is replaced by a higher-dimensional group (or sphere) $\Gamma_n$
so that, outside the fixed-point set, we have a fibration. We
shall refer to such fibrations as {\it branched fibrations}.

Consider the case of the second equation in (\ref{Hopffibration})
involving $U(1)$. The local model here is the action of $U(1)$ on
$C^2 = R^4$ via complex scalars. The quotient is $R^3$ with $S^2$
being the boundary of the branch point. The geometry of the
branched $U(1)$-fibration $R^4 \to R^3$ is fundamental in physics
where it describes the geometry of a magnetic monopole. It was a
major discovery by Dirac that the quantization of electric charge
could be explained (in modern terms) by the $U(1)$-bundle above,
over the complement of the point magnetic source at the origin.
$R^4$ in this situation is now referred to as the Kaluza-Klein
model of the Dirac monopole.

More generally, in current physical theories where space-time is
viewed as having higher dimension than $4$, a $U(1)$-action with a
fixed manifold of codimension $4$ is viewed as providing a charge
on the branch locus (which has codimension $3$ in the quotient).
Examples of such situations were, for example, studied in detail in
\cite{AtWi} and provided some of the early motivation for this
paper.

In a similar way an $Sp(1)$-action with a fixed-point set of
codimension $8$ (and the standard action on $H^2 = R^8$)
gives a branched fibration carrying an $Sp(1) =
SU(2)$-``charge'' on the branch locus, which has codimension $5$
in the quotient.

{Finally} the last equation in (\ref{Hopffibration}) gives a similar
story for branched fibrations with fibre $S^7$ and branch locus
having codimension $9$ in the quotient. Notice that, in this case,
the fibration is not a group action.

In all these cases the quotient manifold has the induced topology
but {\it not} the induced differentiable structure. In other
words, it is not true that a differentiable function above, which is
invariant under the group action, is a differentiable function
below. For example for the double cover $w = z^2$ the function
$x^2$, where $x = \Re(z)$, is not a differentiable function of
$\Re(w),\Im(w)$. However, there {\it is} a natural differentiable
structure on the quotient and we shall always use this and refer
to it as the quotient structure. Note that, for holomorphic
functions the invariants are indeed the functions of $w$ and so
the induced holomorphic structure on the quotient agrees with our
differentiable quotient structure.

In the examples of branched fibrations which we shall study there
will be a further group action in addition to the actions of the
type in (\ref{Hopffibration}). These will be actions of
``cohomogeneity one'', i.e.\ an action of a connected
Lie group $G$ on a connected smooth manifold $M$
whose generic orbit has codimension one.
If $G$ and $M$ are compact such an action has a
simple global structure: either there are no
exceptional orbits and we have a fibration over the circle, or
else there are just two exceptional orbits and the quotient is the
closed unit interval, see e.g.\ \cite{Mos}.
In this second case the exceptional orbits have
isotropy groups $K_1$ and $K_2$ and the generic orbit has isotropy
group $K \subset K_1 \cap K_2$
when we consider the isotropy groups
along a suitable path which connects the two
exceptional orbits. Moreover, the homogeneous spaces
$$
K_1/K = S^p\ \ \ {\rm and}\ \ \ K_2/K = S^q
$$
must both be spheres, where $p+1$ and $q+1$ are the codimensions
of the two exceptional orbits. The normal sphere bundles of these
two orbits are the maps
$$
G/K \to G/K_1\ \ \ {\rm and}\ \ \ G/K \to G/K_2\ .
$$
Finally the manifold $M$, with its $G$-action, is entirely
determined by the conjugacy class of the triple of subgroups
$K_1$, $K_2$ and $K$.

The prototype example, which will be analyzed carefully in Section
4, is when $M = CP^2$ and $G = SO(3)$ acting with two exceptional
orbits of codimension $2$, namely $RP^2$ and $S^2$ (the conic
$z_1^2 + z_2^2 + z_3^2 = 0$).
Here $K_1 = S(O(1) \times O(2)) \cong O(2)$ and
$K_2 = SO(2) \times SO(1) \cong SO(2)$
are embedded in $SO(3)$ so that
$K = K_1 \cap K_2 = S(O(1) \times O(1)) \times SO(1) \cong Z_2$
and $p = q = 1$.

\section{Projective Planes}

In addition to the projective lines over the division
algebras $A_n$ we can consider higher-dimensional projective
spaces. For $n=0,1,2$, when $A_n$ is associative, this gives us the
classical projective spaces
$$
RP^m\ ,\ CP^m\ ,\ HP^m\ (m \geq 2)\ .
$$
For $n = 3$ however, $A_3 = O$ (the octonions) is not
associative. In this case it is possible to construct a projective
plane $OP^2$ (the Cayley plane), but not the projective spaces of
higher dimension. In fact the non-associativity of $O$ is related
to the {\it non-Desarguesian} property of $OP^2$, and it is known
(see e.g.\ \cite{VeYo}) that projective spaces of dimension $\geq 3$ must be
Desarguesian.

Just as with projective lines there are two ways of defining a
projective plane over $A_n$. The first is to introduce the affine
plane in the obvious way as pairs $(x,y)$ of points $x$ and $y$
in $A_n$ and then to
compactify this by adding a projective line at infinity. This
defines $P_n$ as a manifold and the Hopf fibration appears
naturally as the fibration of a spherical neighbourhood of the
line at infinity. The fact that two lines meet in one point gets
translated in this way to an assertion about the topology of the
Hopf fibrations, namely that the Hopf invariant (or linking number
of two fibres) is one. This is the fact which is used to show that
the dimension of a division algebra over $R$ must be $2^n$, $n=0,1,2,3$: see
\cite{AdAt} for a short proof.

This construction of $P_n$ does not exhibit its homogeneity, and
this is where an alternative construction is useful. For $n=0,1,2$
the classical approach is to use the $3$-dimensional vector space
over $A_n$ and to define $P_n$ as the space of one-dimensional
subspaces. This gives $P_n$ as a homogeneous space of the relevant
classical group
\begin{equation}
SL(3,R)\ ,\ SL(3,C)\ ,\ SL(3,H)
\label{projectivegroup}
\end{equation}
or of its compact form
\begin{equation}
SO(3)\ ,\ SU(3)\ ,\ Sp(3)\ .
\label{isometrygroup}
\end{equation}
The linear groups consist of projectivities, i.e.\ transformations
preserving lines. The compact groups consist of isometries, where
the projective plane is equipped
with the Riemannian metric which is induced in the natural way from the
Killing form of the group.
For $CP^2$ the full isometry group is the extension of $SU(3)$ by
$Z_2$ of complex conjugation. For the other two cases the isometry
group is connected. In all cases the centre acts trivially so that
it is really the adjoint group that acts effectively. The isotropy
group for the actions of the three groups in
(\ref{isometrygroup}) are
\begin{equation*}
O(2) \cong S(O(1) \times O(2))\ ,\ U(2) \cong S(U(1) \times U(2))\ ,
\ Sp(1) \times Sp(2)\ .
\label{isotropygroup}
\end{equation*}

{For} $n=3$ we cannot use this approach to construct the Cayley
plane as there is no group $SL(3,O)$. However there is a substitute,
both for the linear group and for its compact form, which plays
the part of the fourth term of the sequences
(\ref{projectivegroup}) and (\ref{isometrygroup}). For
(\ref{isometrygroup}) we have the exceptional compact Lie group $F_4$ and
for (\ref{projectivegroup}) we have the non-compact real form $E_6^{-26}$
with character $-26$ of the
exceptional complex Lie group $E_6(C)$. The group of projective
transformations of the Cayley plane has been explicitly determined
by Freudenthal in \cite{Fre2}.

At this point it is easy to discuss the homogeneity of the Hopf
fibrations in (\ref{Hopffibration}). If we fix a point $o$ in
$P_n$ then the isotropy group at $o$ of the connected isometry group of
$P_n$ acts transitively on the dual projective line $o^*$ in $P_n$
(the set of
all antipodal points of $o$ in $P_n$) and on the metric sphere
bundle over $o^*$ in $P_n$ of sufficiently small radius. The
projection of this sphere bundle from $o$ onto $o^*$ is just the
Hopf fibration associated with $A_n$. The
isotropy group of this action at a point in $o^*$ acts
transitively on the corresponding fibre.

The best way to unify all the projective planes $P_n$, and the
associated groups of symmetries, is to introduce Jordan algebras.
For a quite self-contained treatment in the octonionic case
we refer to \cite{Fre2} and \cite{Mur}, the other cases work
analogously and are easier
to deal with. We summarize here the basic facts.

{For} each division algebra $A_n$ we consider the real vector space
$H_n$ of $3 \times 3$ Hermitian matrices over $A_n$. Recall that
in $A_n$ we have a notion of conjugate $x \mapsto \bar{x}$ which
fixes the ``real'' part and changes the sign of the ``imaginary''
part. Note that conjugation is an anti-involution of the algebra.
Then, as usual, a matrix is Hermitian if its conjugate is
equal to its transpose,
$$
\bar{x}_{ij} = x_{ji}\ \ (i,j=1,2,3)\ .
$$
We make $H_n$ into a commutative (but non-associative) algebra by
defining a multiplication
\begin{equation}
X \circ Y = \frac{1}{2}(XY + YX)\ ,
\label{Jordanmult}
\end{equation}
where $XY$ and $YX$ denote usual matrix multiplication.
Together with this multiplication $H_n$ becomes a real Jordan
algebra which we denote by $J_n$.
The unit matrix $I$ acts as an identity.

{For} $n=0,1,2$ the groups in (\ref{isometrygroup}) act on $H_n$ by
$$
X \mapsto AXA^*\ \ \ (A^* = \bar{A}^t)
$$
and preserve the multiplication (\ref{Jordanmult}). Hence they are
automorphisms of the Jordan algebra $J_n$. It can be shown that
modulo their centres they are the full group $\Aut(J_n)$
of automorphisms of $J_n$, except
for $n=1$ when we get the identity component. For $n=3$
the automorphism group $\Aut(J_3)$ provides an explicit model
of the exceptional compact Lie group $F_4$.

{For} $n=0,1,2$ we have a natural embedding of the projective plane
$P_n$ in $H_n$. We just associate to a one-dimensional subspace
of the $3$-dimensional vector space $A_n^3$ over $A_n$
the Hermitian $3 \times 3$ matrix which represents
orthogonal projection onto it. In terms of homogeneous coordinates
$(x_1,x_2,x_3)$, normalized so that
${\parallel} x_1 {\parallel}^2 + {\parallel} x_2 {\parallel}^2
+ {\parallel} x_3 {\parallel}^2 =
1$,
this is given by
$$
x \mapsto X = (x_i\bar{x}_j)_{i,j=1,2,3}\ .
$$
Note that $x$ and
$x\lambda$ with $\lambda$ in $A_n$,
${\parallel} \lambda {\parallel} = 1$, give the same matrix,
and that $X$ satisfies\footnote{The definition of the determinant
for $n=2$ (the quaternionic case) requires a little care, see
Section 7.}
\begin{equation}
\begin{array}{rcl}
\tr X & = & 1 \ ,\\
{\parallel} X {\parallel}^2 & = & 1\ ,\\
\det X & = & 0\ .
\end{array}
\label{propertyX}
\end{equation}
Clearly the image of $P_n$ in $H_n$ is just the orbit of the
diagonal matrix $\Diag(1,0,0)$ under the isometry group\footnote{This
still works for $P_1 = CP^2$, where the isometry group is
disconnected.}.

{For} $n=3$ it can be shown that the
orbit of $\Diag(1,0,0) \in J_3$ under the action of $F_4 = \Aut(J_3)$
provides a model for the Cayley plane $P_3$.
The isotropy group is isomorphic to $Spin(9)$ and hence
$OP^2 = P_3 = F_4/Spin(9)$ as a homogeneous space \cite{Bor}. If
we consider, as discussed above,
the Hopf fibration $S^{15} \to S^8$ as a projection of a
metric sphere bundle over $o^*$ from a point $o$ in $P_3 = OP^2$ onto the
dual projective line $o^*$, and if $Spin(9)$ denotes the isotropy
group of $F_4$ at $o$, then $Spin(9)$ acts transitively on $S^{15}$
with isotropy group $Spin(7)$ and transitively on $S^8$ with
isotropy group $Spin(8)$. Moreover, $Spin(8)$ acts transitively on
the corresponding fibre $S^7$.

In all four cases $J_n$ has three invariant
polynomials of degrees $1,2,3$ as in (\ref{propertyX}), and we
shall use the same notation. This follows from the fact that every
element in $J_n$ can be reduced by $\Aut(J_n)$ to real diagonal
form. If the diagonal entries (the ``eigenvalues'')
are $\lambda_1,\lambda_2,\lambda_3$ then the three invariant
polynomials are
\begin{equation}
\begin{array}{rcl}
\tr & = & \lambda_1 + \lambda_2 + \lambda_3 \ ,\\
{\parallel} \ \  {\parallel}^2 & = & \lambda_1^2 + \lambda_2^2 + \lambda_3^2\ ,\\
\det & = & \lambda_1 \lambda_2 \lambda_3\ .
\end{array}
\label{propertypoly}
\end{equation}
In all cases only the cubic polynomial $\det$ is
invariant under the group of projectivities of
the projective plane $P_n$. The fact that it is actually given by
a polynomial needs to be proved (see Section 7).

\section{Orbit Structures}

As mentioned at the end of Section 2 the manifolds we are
interested in have cohomogeneity one group actions compatible with
the maps we need to construct to prove the identifications
(\ref{quotient}). We proceed to spell these out in detail
beginning with the basic example of $CP^2$, which will be a model for
the others.

We consider the action of $SO(3)$ on $CP^2$ via the natural
embedding $SO(3) \subset SU(3)$. There are two special orbits,
namely $RP^2$ with isotropy group $K_1 = S(O(1) \times O(2)) \cong O(2)$, and
a $2$-sphere $S^2$ (the conic $z_1^2 + z_2^2 + z_3^2 = 0$)
with isotropy group $K_2 =
SO(2) \times 1 \cong SO(2)$. The intersection
\begin{equation}
K = K_1 \cap K_2 = S(O(1) \times O(1)) \times 1 \cong Z_2
\label{isotropyplaneC}
\end{equation}
consists of the two diagonal matrices $\Diag(\lambda,\lambda,1)$
with $\lambda = \pm 1$. Each of the homogeneous spaces $K_1/K$ and
$K_2/K$ is a circle. The generic orbit $SO(3)/K$ is $3$-dimensional
and it fibres over each of the two special orbits $RP^2$ and $S^2$ with
$S^1$ as fibre.

One way to establish the identity
$$
CP^2/O(1) = S^4
$$
of (\ref{quotient}) is to analyze the $SO(3)$-orbit structure of
$S^4$ and compare it with that of $CP^2$. Here we can view $S^4$
as the unit sphere in the vector space of symmetric $3 \times 3$ real
matrices of trace zero equipped with its usual norm.
The orbits are determined by the three
real eigenvalues $\lambda_1 \leq \lambda_2 \leq \lambda_3$ with
$\lambda_1 + \lambda_2 + \lambda_3 = 0$. There are two special
orbits $\lambda_1 = \lambda_2$ and $\lambda_2 = \lambda_3$ each
of which is an $RP^2$, while the generic orbit is the real flag
manifold of all full flags in $R^3$.
Thus the three isotropy groups are
\begin{equation}
K_1^\prime = S(O(1) \times O(2))\ ,\ K_2^\prime = S(O(2) \times O(1))\
,\ K^\prime = K_1^\prime \cap K_2^\prime = S(O(1)^3)\ .
\label{isotropysphere4}
\end{equation}

Since the $4$-manifold with its $SO(3)$-orbit structure is determined
by the conjugacy class of the triples of isotropy groups,
comparison of (\ref{isotropyplaneC}) and (\ref{isotropysphere4})
shows that there is a natural map
$$
CP^2 \to S^4
$$
compatible with the two $SO(3)$-actions. Moreover this map identifies
the special orbit $RP^2$ in $CP^2$ with one of the two $RP^2$ in
$S^4$, while outside this we have a double covering, given by the action
of $O(1)$ on $CP^2$ as complex conjugation.

The orbit structures of $HP^2$ and $OP^2$ are quite similar as are
those of the corresponding spheres. We consider first the action
of $SU(3)$ on $HP^2$ via the embedding $SU(3) \subset Sp(3)$, and on
$S^7$ as the unit sphere in the vector space of Hermitian $3 \times 3$
complex matrices of trace zero equipped with its usual norm.
For $S^7$ we find two copies of $CP^2$ as
special orbits and the complex flag manifold of all full flags
in $C^3$ as generic orbit, so
that the isotropy groups are
\begin{equation}
K_1^\prime = S(U(1) \times U(2))\ ,\ K_2^\prime = S(U(2) \times U(1))\
,\ K^\prime = K_1^\prime \cap K_2^\prime = S(U(1)^3)\ .
\label{isotropysphere7}
\end{equation}
For $HP^2$ the two special orbits are $CP^2$ and a circle bundle
over the dual $CP^2$ which is a $5$-dimensional sphere $S^5$
(as one sees by using all quaternion lines,
i.e.\ $4$-spheres, determined by complex lines). The isotropy groups
are
\begin{equation*}
K_1 = S(U(1) \times U(2))\ ,\ K_2 = SU(2) \times 1\ ,\ K = K_1
\cap K_2 = S(U(1)^2) \times 1\ .
\label{isotropyplaneH}
\end{equation*}
Comparison with (\ref{isotropysphere7}) shows the existence of a
map\footnote{The proof outlined here is essentially that of
\cite{Arn}, \cite{AtWi}.}
$$
HP^2 \to S^7
$$
compatible with the two $SU(3)$-actions. It identifies the $CP^2$
in $HP^2$ with one of the two $CP^2$ in $S^7$, and on the
complement $HP^2 \setminus CP^2$ it is an $S^1$-bundle over the
complement $S^7 \setminus CP^2$. The $SU(3)$ determines a
maximal subgroup $U(3)$ of $Sp(3)$, and the central $U(1)$ in this
$U(3)$ acts trivially on the $CP^2$ and gives the fibres on
$HP^2 \setminus CP^2$.

{Finally} consider the action of $Sp(3)$ on the Cayley plane $OP^2$
and on $S^{13}$, the unit sphere in the space of Hermitian $3
\times 3$ quaternion matrices with trace zero equipped with its usual norm.
Again we have two
special orbits in $S^{13}$, both copies of $HP^2$, and the generic
orbit is the quaternionic flag manifold of all full flags in
$H^3$, so that the
isotropy groups are
\begin{equation}
K_1^\prime = Sp(1) \times Sp(2)\ ,\ K_2^\prime = Sp(2) \times Sp(1)\
,\ K^\prime = K_1^\prime \cap K_2^\prime = Sp(1)^3\ .
\label{isotropysphere13}
\end{equation}
For the Cayley plane, $HP^2$ is clearly a special orbit. The
generic orbit is just the normal sphere bundle of $HP^2$ in $OP^2$
with fibre $S^7$. This is the unit sphere in the normal $R^8$
which is a representation of the isotropy group $K_1 = Sp(1)
\times Sp(2)$. Note that this representation is {\it not} the
standard action on $H^2$ since $Sp(1)$ acts trivially and only
$Sp(2)$ acts in the standard manner, so that the generic isotropy
group is $K = Sp(1) \times Sp(1) \times 1$. By considering all the
Cayley lines ($8$-spheres) determined by quaternion lines we see
that the other special orbit is fibred over the dual $HP^2$ with
$3$-sphere fibres and therefore is an $11$-dimensional sphere
$S^{11}$. Hence the isotropy groups are
\begin{equation*}
K_1 = Sp(1) \times Sp(2)\ ,\ K_2 = Sp(2) \times 1\ ,\ K = K_1
\cap K_2 = Sp(1)^2 \times 1\ .
\label{isotropyplaneO}
\end{equation*}
Comparison with (\ref{isotropysphere13}) then shows the existence
of a map $$ OP^2 \to S^{13} $$ compatible with the action of
$Sp(3)$. $HP^2$ is the branch locus and outside this we have a
$3$-sphere fibration. The $Sp(3)$ is contained in a maximal
subgroup $Sp(3)Sp(1)$ in $F_4$, where the $Sp(1)$ centralizes the
$Sp(3)$ in $F_4$, see e.g.\ \cite{BoSi}. This $Sp(1)$ fixes the
$HP^2$ in $OP^2$ pointwise, and the orbits through the other
points in $OP^2$ are just the $3$-spheres of the fibration. This
establishes the identity $OP^2/Sp(1) = S^{13}$ in
(\ref{quotient}).

We have thus established (\ref{quotient}) from a complete
description of the relevant orbit structures.
In the next section
we shall formulate and prove Theorem A, a more explicit version of
(\ref{quotient}) which does not rely on such a detailed knowledge
of the orbit structure, and which provides such an explicit map.

The orbit structure of these group actions on the projective
planes have also been studied in \cite{PuRi} in relation to
homotopy theory.

\section{An Explicit Map}

In Section 3 we saw that in all four cases the projective plane
$P_n$ has a natural embedding in $H_n$, the vector space of $3
\times 3$ Hermitian matrices over the division algebra $A_n$.
Moreover the formulae (\ref{propertyX}) show that $P_n$ lies in
the hyperplane $H_n(1)$ given by $\tr X = 1$ and on the sphere
${\parallel}X{\parallel}^2 = 1$. Note that the intersection of the
hyperplane and the sphere is the sphere of one lower dimension
with centre $I/3$ and radius $\rho$ where $I$ is the unit matrix
$\Diag(1,1,1)$ and $\rho^2 = 2/3$. Since
$$
\dim H_n = 3(2^n+1)
$$
we thus have an embedding
\begin{equation}
P_n \subset S^{d(n)}\ ,\ d(n) = 3\cdot 2^n + 1\ .
\label{planeinsphere}
\end{equation}
For $n=0,1,2$ these are the classical embeddings referred to in
(\ref{branchloci}) and for $n=3$ we have the corresponding one for
the Cayley plane.

{For} $n\geq 1$ we have a natural inclusion of algebras
$$
A_{n-1} \subset A_n
$$
and hence, using the Euclidean metric given by the norm, an
orthogonal projection $A_n \to A_{n-1}$ which induces a projection
\begin{equation*}
\pi_n : H_n \to H_{n-1}\ .
\label{projection}
\end{equation*}
Note that $\pi_n$ commutes with the trace and hence it maps
$P_n$ into the affine hyperplane $H_{n-1}(1)$ in $H_{n-1}$. In
this hyperplane we will choose the point $I/3$ as centre and use
the shifted coordinate
$$
\tilde{X} = X - I/3
$$
for $H_{n-1}(0)$, the linear subspace of $H_{n-1}$ given by
$\tr X = 0$. Restricting $\pi_n$ to $P_n \subset H_n$ we get a map
$$
\pi_n : P_n \to H_{n-1}(1)
$$
and a shifted map $\tilde{\pi}_n : P_n \to H_{n-1}(0)$. The following lemma will be
crucial for our construction.

\medskip
{\sc Lemma 1.} {\it The matrix $I/3$ does not lie in the image $\pi_n(P_n)$, so that
$\tilde{\pi}_n(X) \neq 0$ for all $X \in P_n$.}

\medskip
We postpone the proof till later. Assuming that this lemma is true we can
rescale the maps $\tilde{\pi}_n$ to define a map
\begin{equation}
f_n : P_n \to S^{d(n-1)}
\label{mapf}
\end{equation}
given by
\begin{equation}
f_n(X) = \frac{1}{3}I +
\rho\frac{\tilde{\pi}_n(X)}{{\parallel}\tilde{\pi}_n(X){\parallel}}\
. \label{mapfX}
\end{equation}
Thus $f_n(P_n)$ lies in the sphere $S^{d(n-1)}$ of radius $\rho$ in
$H_{n-1}(1)$ centred at $I/3$.

The map (\ref{mapf}) will be the map inducing the diffeomorphism
of (\ref{quotient}), but before formulating Theorem A we shall
need a few further properties of $f_n$. First we note that, when
$X \in P_{n-1} \subset P_n$, $\pi_n(X) = X$ and
${\parallel}X{\parallel} = 1$, so that
${\parallel}\tilde{\pi}_n(X){\parallel} = \rho$ and $f_n(X) = X$.
Hence, when restricted to $P_{n-1}$, the map $f_n$ is just the
standard embedding (\ref{planeinsphere}) for $n-1$.

Let us denote by $G_n$ the groups of isometries of $P_n$ which we
already discussed above, namely
\begin{equation}
\begin{array}{rcl}
G_0 & = & SO(3)\ ,\\
G_1 & = & SU(3)\ ,\\
G_2 & = & Sp(3)\ ,\\
G_3 & = & F_4\ .
\end{array}
\label{groupG}
\end{equation}
The space $H_n$ of Hermitian matrices is a representation of
$G_n$, which splits off a trivial factor
(corresponding to the trace).
When restricted to $G_{n-1}$ ($n=1,2,3$) it decomposes as
\begin{equation}
H_n = H_{n-1} \oplus H_{n-1}^{\perp}\ ,
\label{Hdecomposed}
\end{equation}
so that the projection $\pi_n : H_n \to H_{n-1}$ is compatible
with the action of $G_{n-1}$. Thus, assuming Lemma 1, the map
$f_n$ of (\ref{mapf}) also commutes with the action of $G_{n-1}$.

In addition we have an action of the group $\Gamma_{n-1}$ of
elements of norm one
of the division algebra $A_{n-1}$ on $H_n$ which commutes with the
action of
$G_{n-1}$ and preserves the fibres of $\pi_n$. Recall that $G_n$
is, modulo its centre, the automorphism group $\Aut(J_n)$ of the Jordan
algebra $J_n$ except for $n=1$ when $\Aut(J_1)$ has another connected
component induced by complex conjugation. Note that this complex
conjugation does not extend to an automorphism of $J_2$ but to an
anti-automorphism. Then $\Gamma_{n-1}$ can be viewed as the
centralizer of the connected component of $\Aut(J_{n-1})$ in
$\Aut(J_n)$. Explicitly we have
\begin{itemize}
\item[$n=1$:] $\Gamma_0 = O(1) \cong Z_2$ acts by complex conjugation on
$H_1$;
\item[$n=2$:] $\Gamma_1 = U(1)$ acts by conjugation on
$H_2$ with respect to the elements of norm one of $C \subset H$;
\item[$n=3$:] $\Gamma_2 = Sp(1)$. The action of $Sp(1)$ on $H_3$
cannot be described in a similar fashion because of the
non-associativity of the octonions, but the construction of $Sp(1)$
as a subgroup of $\Aut(J_3) = F_4$ is quite simple. Consider a root space
decomposition of the Lie algebra of $F_4$ such that the Lie algebra
of $G_2 = Sp(3)$ is determined by the two short simple roots and
the adjacent long simple root. Then the maximal root determines
the Lie algebra of $\Gamma_2 = Sp(1)$.
\end{itemize}

\noindent Note that $\Gamma_{n-1}$ acts trivially on the summand $H_{n-1}$
in the decomposition (\ref{Hdecomposed}) and hence trivially on
the projective plane $P_{n-1} \subset H_{n-1}$.
The action of $\Gamma_0 = O(1)$ on $H_0^{\perp} \cong R^3$ is just
by scalar multiplication by $\pm 1$, the one of $\Gamma_1 =
U(1)$ on $H_1^{\perp} \cong C^3$ is the scalar action
$(\lambda,z) \mapsto \lambda^2z$ ($\lambda \in U(1) \cong S^1
\subset C$, $z \in C^3$), and the one of $\Gamma_2 = Sp(1)$ on
$H_2^{\perp} \cong H^3$ is by right multiplication.
The action of $\Gamma_{n-1}$ on the normal bundle of $P_{n-1}$ in
$P_n$ is obtained from the one on $H_{n-1}^{\perp}$ by suitable restriction.

{Finally} we are in a position to formulate our promised refinement
of (\ref{quotient}).

\medskip
{\sc Theorem A.} {\it For $n = 1,2,3$ the map
$$
f_n : P_n \to S^{d(n-1)}
$$
defined by (\ref{mapfX}) induces a diffeomorphism
$$
P_n/\Gamma_{n-1} \approx S^{d(n-1)}\ .
$$
Moreover, this diffeomorphism is compatible with the natural
action on $P_n/\Gamma_{n-1}$ and $S^{d(n-1)}$
of the group $G_{n-1}$ of (\ref{groupG}).}

\medskip
{\sc Remark.}
We have already observed that the action of
$\Gamma_{n-1}$ on the normal bundle of $P_{n-1}$ in $P_n$ is just
the scalar action of the relevant field
(for $n=2$ it is the square of the scalar action),
so that the quotient
$P_n/\Gamma_{n-1}$ is indeed a manifold, branched along $P_{n-1}$
in the sense of Section 2. Outside the branch locus the action of
$\Gamma_{n-1}$ is free (except for $n=2$ where each orbit is a
double covering of the circle, in which case we may consider
$\Gamma_1/Z_2 \cong \Gamma_1$ to get a free action).

\smallskip
The case $n=1$ of Theorem A gives the known diffeomorphism of $P_1
= CP^2$ modulo complex conjugation with $S^4$. Since this case is
basic to the others we shall prove it first.

Consider, as a preliminary, the complex projective line $CP^1$
viewed as embedded by idempotents in the affine space $R^3$ of
Hermitian $2 \times 2$ complex matrices of trace one. The projection onto
the real symmetric matrices of trace one is easily seen to be the
standard projection of $S^2$ onto the disc $D^2$ (with centre
$I/2$ and radius $\sigma$ with $\sigma^2 = 1/2$), which identifies
conjugate points. Note that the $O(2)$-orbits on $S^2$, the ``circles of
latitude'', go into the concentric circles in the disc.

We are now ready to look at $CP^2$ and the projection $\pi_1$. For
every $RP^1$ in $RP^2$ its complexification is a $CP^1$ in $CP^2$
and these fill out $CP^2$, intersecting only at points of $RP^2$.
Thus $CP^2 \setminus RP^2$ is fibred over the dual $RP^2$ with
fibre $CP^1 \setminus RP^1 = S^2 \setminus S^1$ (two open discs).
Alternatively it is fibred over $S^2$ with fibre an open disc.
Moreover, $SO(3)$ acts transitively on $RP^2$ (or $S^2$), the base
of this fibration, and the isotropy group
$S(O(1) \times O(2)) \cong O(2)$ acts on each
fibre. Because the projection $\pi_1$ is compatible with this
action of $SO(3)$ it is entirely determined by its restriction to
a single fibre. But such a fibre can be taken to be given
by the equation $z_3 = 0$, and so we are reduced to studying the
projection of $CP^1$ which we have just done above. The first
implication of this is to establish Lemma 1 for $n=1$, because
the scalar $3 \times 3$ matrix $I/3$ cannot lie in the subspace of
$2 \times 2$ matrices (with zeroes in the third row and column).
Moreover, the orbit analysis of the $2$-dimensional case shows that
the map $f_1$ sends the $SO(3)$-orbits of $CP^2$ (modulo
conjugation) diffeomorphically onto the $SO(3)$-orbits of $S^4$,
thus proving Theorem A for $n=1$.

To deal with the cases $n=2,3$ of Theorem A, i.e.\ with $HP^2$ and
$OP^2$, we shall choose appropriate embeddings of $CP^2$ in the
higher projective planes and then use Theorem A for $CP^2$. The
embeddings we want are not the standard ones given by the original
inclusions $C \subset H \subset O$, but are suitable conjugates of
these. Thus, for $H$, we choose the embedding $C \to H$ which
takes $i \in C$ into $j \in H$, which gives a second embedding of
$CP^2$ into $HP^2$, which we will denote by $CP^2(j)$ to
distinguish it from the original $CP^2$. Note that
\begin{equation}
CP^2(j) \cap CP^2 = RP^2\ .
\label{intersectH}
\end{equation}
Similarly we will choose a third embedding of $CP^2$ into $OP^2$
(not coming from $HP^2$). Consider the sphere $S^6$ of imaginary
elements of norm one in $O = R^8$. This contains the $S^2$ of imaginary
quaternions of norm one. Choose an element $e \in S^6$ of norm one
which lies in $R^4$
orthogonal to $H \subset O$, and embed $C$ in $O$ by sending $i$
to $e$. Since the exceptional compact Lie group\footnote{Note that this is
not to be confused with the group $G_2 = Sp(3)$ of the sequence in
(\ref{groupG}).} $G_2$ of automorphisms of the octonions acts
transitively on $S^6$ we can find $g \in G_2$ which takes $i$ into
$e$. Since $G_2$ acts naturally on the exceptional Jordan algebra
$J_3$ (so that $G_2 \subset F_4$) we get a copy $g(CP^2) \subset
OP^2$. Note that
\begin{equation}
g(CP^2) \cap CP^2 = RP^2\ .
\label{intersectO}
\end{equation}

Because our elements $j$ and $e$ were chosen orthogonal to $i$ it
follows that the projections $\pi$ of $HP^2$ and $OP^2$ restrict
to the standard projection of $CP^2(j)$ and $g(CP^2)$. Because of
(\ref{intersectH}) and (\ref{intersectO}) these copies of $CP^2$
are transversal to the orbits of the relevant groups $SU(3)$ and
$Sp(3)$. More precisely, each orbit of the larger group
intersects our $CP^2$ in just one $SO(3)$-orbit. This follows by
examining the corresponding groups. Consider first the embedding
$CP^2(j) \subset HP^2$. This is an orbit of a copy of $SU(3)$
which we may denote by $SU(3)_j$. Clearly this intersects the
original $SU(3) \subset Sp(3)$ precisely in $SO(3)$. Similarly
$g(CP^2) \subset OP^2$ is an orbit of $g(SU(3))$ and we need to
check that
$$
SU(3) \cap g(SU(3)) = SO(3)\ .
$$
But this is clear because the intersection must preserve
$$
CP^2 \cap g(CP^2) = RP^2\ .
$$

This correspondence between the $SO(3)$-orbits on $CP^2$ and the
orbits of $SU(3)$ on $HP^2$ and of $Sp(3)$ on $OP^2$ shows that
Lemma 1 for $n=2,3$ follows from the case $n=1$, which we have
already proved. The correspondence also shows that the map $f_n$
induces a diffeomorphism on the one-dimensional space (interval)
of orbits. But each orbit in $P_n$ is known and outside the branch
locus it is just fibred over the corresponding fibre in
$S^{d(n-1)}$. Together with the local behaviour near the branch
locus this completes the proof of Theorem A.

\section{The Complexified Projective Planes}

As mentioned in Section 1 the four projective planes $P_n$
have natural complexifications $P_n(C)$ as complex
projective algebraic varieties.  These have a ``real structure'',
i.e.\ an anti-holomorphic involution $\tau$, which has $P_n$ as
the real part (fixed by $\tau$). We shall now examine these
varieties in greater detail and study their symmetries.

Recall that $P_n \subset H_n$, the space of Hermitian $3 \times 3$
matrices over the division algebra $A_n$. It is the orbit of the
diagonal matrix $\Diag(1,0,0)$ under the compact Lie group $G_n$
(listed in (\ref{groupG})). Note that $P_n \subset H_n(1)$, the
affine subspace of matrices of trace one. If we denote by $\cP_n$
the real projective space of the vector space $H_n$, so that
$\dim \cP_n = 3\cdot 2^n +2$, then we can identify $\cP_n$ with
the projective completion of the real affine space $H_n(1)$ and
we have $P_n \subset \cP_n$.

The group $G_n$ of isometries of $P_n$ acts on the affine space
$H_n(1)$. This action extends to an action of a larger group
$\cG_n$ on the projective space $\cP_n$ which preserves $P_n$.
This induces on $P_n$ its group of projectivities. For $n=0,1,2$
we have
$$
\cG_n = SL(3,A_n)
$$
as noted in (\ref{projectivegroup}), while $\cG_3$ is the
non-compact real form $E_6^{-26}$
of the exceptional complex Lie group $E_6(C)$.
In all cases $\cG_n$ has a natural irreducible
representation on the real vector space $H_n$, which splits off a
trivial factor (corresponding to the trace) when restricted to the
compact subgroup $G_n$.

{For} $n=0$ the representation of $\cG_0 = SL(3,R)$ in $H_0 = R^6$
is via the symmetric square $S^2(R^3)$, and the embedding $P_0
\subset \cP_0$ is the embedding
$$
RP^2 \subset RP^5
$$
induced by the diagonal (squaring) map $R^3 \to S^2(R^3)$ given by
$x \mapsto x^2$. It is the (real) Veronese embedding, and we shall
use the same term for all $n$.

Note that, if $V = R^3$, then $SL(V)$ has two inequivalent irreducible
representations of dimension $6$, namely $S^2(V)$ and $S^2(V^*)$,
where $V^*$ denotes the dual vector space of $V$.
These become equivalent when restricted to
$SO(3)$. A similar story holds for all $n$, so that these Veronese
embeddings come in dual pairs
\begin{equation*}
P_n \subset \cP_n\ ,\ P_n^* \subset \cP_n^*\ ,
%\label{dualVeronese}
\end{equation*}
where $P_n^*$ is the dual projective plane of $P_n$, representing
its projective lines. Reduction to the compact group gives an
identification
$$
P_n \to P_n^*
$$
by associating to each point $p$ of the plane $P_n$ the
``opposite'' projective line, which may be defined as the set of
all antipodal points $q$ on a closed geodesic through $p$,
the ``polar'' of $p$ in $P_n$.

We are now in a position to complexify everything. We get a
complex Lie group $\cG_n(C)$ acting on the complex projective
space $\cP_n(C)$. Since $P_n$ is an orbit of $\cG_n$ we get as
complexification an orbit $P_n(C)$ of $\cG_n(C)$, which defines the
complexified projective plane. It is a complex projective manifold
with
$$
\dim_C P_n(C) = \dim_R P_n = 2^n\ .
$$

Since $\cP_n(C)$ has a natural Hermitian metric we can define the
maximal compact subgroup $\hat{G}_n$ of $\cG_n(C)$. This
preserves the induced K\"{a}hler metric on $P_n(C)$ and it clearly
contains $G_n$. Note that $P_n(C)$ equipped with this induced
K\"ahler metric is a Hermitian symmetric space.

Explicitly, the groups $\cG_n(C)$ and $\hat{G}_n$ are given by
the following table, where for greater clarity we also list the
groups $G_n$ and $\cG_n$:
\begin{equation}
\renewcommand{\arraystretch}{1.5}
\begin{array}{|c|c|c|c|c|}
\hline
n & 0 & 1 & 2 & 3 \\
\hline
\cG_n(C) & SL(3,C) & SL(3,C) \times SL(3,C) & SL(6,C) & E_6(C) \\
\hline
\cG_n & SL(3,R) & SL(3,C) & SL(3,H) & E_6^{-26} \\
\hline
\hat{G}_n & SU(3) & SU(3) \times SU(3) & SU(6) & E_6 \\
\hline
G_n & SO(3) & SU(3) & Sp(3) & F_4 \\
\hline
\end{array}\ .
\label{table}
\end{equation}

The compact complex manifolds $P_n(C)$ are necessarily homogeneous
spaces also of the maximal compact subgroup $\hat{G}_n$ of the
complex Lie group $\cG_n(C)$. Explicitly we have
\begin{equation*}
\renewcommand{\arraystretch}{1.3}
\begin{array}{rclcl}
P_0(C) & = & SU(3)/S(U(1)\times U(2)) & = & CP^2\ ,  \\
P_1(C) & = & SU(3)^2/S(U(1)\times U(2))^2 & = & CP^2 \times CP^2 \ , \\
P_2(C) & = &  SU(6)/S(U(2)\times U(4)) & = & Gr_2(C^6)\ ,  \\
P_3(C) & = & E_6/Spin(10)U(1) \ . & &
\end{array}
%\label{complexplanes}
\end{equation*}
Here $Gr_2(C^6)$ denotes the complex Grassmannian of 2-planes in
$C^6$. The identification of the isotropy groups is easy in the
classical cases ($n=0,1,2$) and follows from the representation
theory of $E_6$ and $F_4$ for the last case \cite{Ada}. The
embeddings $P_n(C) \to \cP_n(C)$ are well-known classical
embeddings for $n=0,1,2$. For $n=0$ it is the (complex) Veronese
embedding $CP^2 \to CP^5$, for $n=1$ it is the Segre embedding
$CP^2 \times CP^2 \to CP^8$, and for $n=2$ it is the Pl\"{u}cker
embedding $Gr_2(C^6) \to CP^{14}$.

The homogeneous space $P_n$ and $P_n(C)$, together with their
isometry groups, are essentially given by the first two rows of
Freudenthal's magic square
\begin{equation}
\renewcommand{\arraystretch}{1.5}
\begin{array}{|c|c|c|c|}
\hline
\cs\co(3) & \cs\cu(3) & \cs\cp(3) & \cf_4 \\
\hline
\cs\cu(3) & \cs\cu(3) \oplus \cs\cu(3) & \cs\cu(6) & \ce_6 \\
\hline
\cs\cp(3) & \cs\cu(6) & \cs\co(12) & \ce_7 \\
\hline
\cf_4 & \ce_6 & \ce_7 & \ce_8 \\
\hline
\end{array}\ ,
\label{magic}
\end{equation}
see for example \cite{Fre1}. According to Freudenthal, the entries in
the first row describe $2$-dimensional elliptic geometries, in the second row
$2$-dimensional projective geometries, in the third row $5$-dimensional
symplectic geometries, and in the last row metasymplectic
geometries.

We now want to look at the action of $G_n$ on $P_n(C)$ and study
its orbit structure. We already know that $P_n$ is one orbit, say
$P_n = G_n/M_n$ where $M_n$ is given by
\begin{equation*}
\begin{array}{lllll}
M_0 & = & S(O(1) \times O(2)) & \cong & O(2)\ ,\\
M_1 & = & S(U(1) \times U(2)) & \cong & U(2)\ ,\\
M_2 & = & Sp(1) \times Sp(2) \ ,\\
M_3 & = & Spin(9)\ .
\end{array}
%\label{isotropyM}
\end{equation*}
Since $P_n(C)$ is the
complexification of $P_n$ the normal bundle $N_n$ is isomorphic to
the tangent bundle and hence the action of $M_n$ on $N_n$ is just
the representation of $M_n$ on the quotient of the Lie algebras
\begin{equation}
L(G_n)/L(M_n)\ .
\label{repofM}
\end{equation}
But for all these representations $M_n$ is transitive on the unit
sphere. It follows that the generic orbit of the action of
$G_n$ on $P_n(C)$ has
codimension one. Moreover the generic isotropy groups are just the
isotropy groups of $M_n$ on the unit sphere in (\ref{repofM}).
Hence the generic orbits are
\begin{equation*}
\begin{array}{rcl}
n = 0 & : & SO(3)/O(1)\ ,\\
n = 1 & : & SU(3)/U(1)\ ,\\
n = 2 & : & Sp(3)/Sp(1) \times Sp(1)\ ,\\
n = 3 & : & F_4/Spin(7)\ .
\end{array}
%\label{genericorbit}
\end{equation*}
For $n = 3$ this follows from the well-known fact that $S^{15} =
Spin(9)/Spin(7)$, where the action of $Spin(9)$ on $S^{15} \subset
R^{16}$ is via its irreducible spin representation \cite{Bor}.

In addition to the special orbit $P_n$ in $P_n(C)$ there must be
another special orbit. Since $G_n$ acts on the affine space $H_n(1)$ it
leaves invariant the hyperplane section of $P_n(C)$ at infinity,
which we already denoted in Section 1 by $P_n(\infty)$.

After these preliminaries on the spaces $P_n(C)$ and the groups
acting on them we now want to explicitly construct the maps
referred to in (\ref{mapphin}). The method is very similar to that
we used to construct the maps $f_n$ of Theorem A, the essential
point being a projection onto a linear space which misses the
origin (Lemma 1) and hence can be normalized to map to the sphere.
For $P_n$ our projection was the orthogonal projection
$$
\pi_n : H_n \to H_{n-1}
$$
restricted to $P_n \subset H_n$. For such a method we need now to
replace the projective embedding $P_n(C) \to \cP_n(C)$ by an
embedding in a Euclidean space. Fortunately such an embedding
exists for every complex projective algebraic variety which is a
homogeneous space of a compact Lie group. We just have to pick an
appropriate orbit in the Lie algebra. Thus we can find embeddings
\begin{equation*}
P_n(C) \to L(\hat{G}_n)
%\label{embedPnC}
\end{equation*}
compatible with the $\hat{G}_n$-action, where $\hat{G}_n$ is the
group of isometries of $P_n(C)$ given by table (\ref{table}).

We can be more precise. Note that $G_n \subset \hat{G}_n$ and that
$\hat{G}_n/G_n$ is the compact dual of the non-compact symmetric
space $\cG_n/G_n$, since $\hat{G}_n$ and $\cG_n$ are both real forms
of the same complex Lie group $\cG_n(C)$. The compact symmetric space
$\hat{G}_n/G_n$ is just the space that parametrizes
all real $P_n$ in $P_n(C)$, whereas
the non-compact symmetric space $\cG_n/G_n$ is
the space that parametrizes the standard metrics on $P_n$ (up to overall
scale) and so can be identified with the open set
$$
H_n^+(1) \subset H_n(1)
$$
consisting of matrices which are positive-definite (i.e.\ all
$X$ in $H_n(1)$ for which
the eigenvalues $\lambda_1,\lambda_2,\lambda_3$ as in
(\ref{propertypoly}) are positive). Here we choose the diagonal
matrix $I/3$ to define our base metric (base point of
$\cG_n/G_n$). It follows that there is a natural isomorphism
between the tangent spaces
\begin{equation}
L(\hat{G}_n)/L(G_n) \cong L(\cG_n)/L(G_n) \cong H_n(0)\ ,
\label{tangentspace}
\end{equation}
where $H_n(0)$ is the vector space obtained from
the affine space $H_n(1)$ with $I/3$ as its origin.

$P_n$ is the $G_n$-orbit of the matrix $\Diag(1,0,0)$ in $H_n(1)$ or
equivalently the $G_n$-orbit of $\Diag(2/3,-1/3,-1/3)$ in $H_n(0)$.
From (\ref{tangentspace}) we get an orthogonal decomposition
\begin{equation}
L(\hat{G}_n) = L(G_n) \oplus H_n(0)\ .
\label{decompofg}
\end{equation}
The $G_n$-orbit $P_n$ in $H_n(0)$ then generates the
$\hat{G}_n$-orbit $P_n(C)$ in $L(\hat{G}_n)$. The projection
$$
\sigma_n : L(\hat{G}_n) \to H_n(0)
$$
defined by (\ref{decompofg}) restricts to give a map
\begin{equation}
\sigma_n : P_n(C) \to H_n(0)\ .
\label{mapsigma}
\end{equation}
Parallel to Lemma 1 we now have

\medskip
{\sc Lemma 2.} {\it The image of the map $\sigma_n$ in
(\ref{mapsigma}) does not contain $0$.}

\medskip
Assuming for the moment the truth of Lemma 2 we can now define a
map
\begin{equation*}
\varphi_n : P_n(C) \to S^{d(n)}\ ,\ d(n) = 3\cdot 2^n + 1\ ,
%\label{mapphi}
\end{equation*}
by the normalization
\begin{equation}
\varphi_n(Z) = \rho\frac{\sigma_n(Z)}{{\parallel} \sigma_n(Z)
{\parallel}}\ ,
\label{mapphiZ}
\end{equation}
where $\rho$ with $\rho^2 = 2/3$ is inserted
to ensure that, on restriction to
$P_n \subset P_n(C)$, the map $\varphi_n$ coincides with the
standard inclusion $P_n \to H_n(0)$.

We are now in a position to formulate the main result of this
section, refining (\ref{mapphin}),

\medskip
{\sc Theorem B.} {\it The map $\varphi_n : P_n(C) \to S^{d(n)}$
defined by (\ref{mapphiZ}) for $n=0,1,2,3$ is a fibration outside
the branch locus $P_n$ and the hypersurface $P_n(\infty)$ which
gets mapped to the antipodal $P_n$. The fibres are the norm one
elements
$\Gamma_n$ of $A_n$ (namely the spheres $S^0,S^1,S^3,S^7$) and
$\varphi_n$ commutes with the action of $G_n$.}

\medskip
{\sc Remark.} The case $n = 0$ of Theorem B coincides with the
case $n=1$ of Theorem A, since $P_0(C) = P_1 = CP^2$ is the
complexification of $RP^2$ and
$$
L(\hat{G}_1) \cong \cs\cu(3) \cong iH_1(0)\ .
$$
Note that Lemma 2 for $n=0$ reduces to Lemma 1 for $n=1$. We can
therefore use the case $n=0$ of Lemma 2 and Theorem B for the
other three cases. We will simply use the natural inclusions
\begin{equation}
P_0(C) \subset P_1(C) \subset P_2(C) \subset P_3(C)
\label{inclusionP}
\end{equation}
and the corresponding inclusions
\begin{equation*}
H_0(0) \subset H_1(0) \subset H_2(0) \subset H_3(0)\ .
%\label{inclusionH}
\end{equation*}
The maps $\sigma_n$ and $\varphi_n$ are compatible with these
inclusions. Moreover, for $n=1,2,3$,
the $G_n$-orbits in $P_n(C)$ intersect
$P_{n-1}(C)$ in the $G_{n-1}$-orbits. Thus each inclusion in
(\ref{inclusionP}) induces naturally a diffeomorphism from the space
of $G_n$-orbits onto the space of $G_{n-1}$-orbits
(which are both closed intervals).
Since the property in Lemma 2 of not containing
$0$ is $G_n$-invariant it is a property of orbits. Lemma 2 for $n
\geq 1$ therefore follows from the case $n = 0$ (since the
$\sigma_n$-image of any $G_n$-orbit of $P_n(C)$ contains the
corresponding $\sigma_0$-image). Moreover, for the same reasons we
see that $\varphi_n$ maps the generic $G_n$-orbit in $P_n(C)$ onto
a generic $G_n$-orbit in $S^{d(n)}$ in a smooth manner (i.e.\ by a
diffeomorphism of the orbit parameter). The identification of the
fibres as $\Gamma_n$ follows from their explicit description.
Finally, the branch locus $P_n$ being preserved by $\varphi_n$,
the other exceptional orbit $P_n(\infty)$ in $P_n(C)$ must go to
the other (dual or antipodal) $P_n$ in $S^{d(n)}$. This completes
the proof of Theorem B.

In the appendix we will give more information about the map
from $P_n(\infty)$ onto this dual $P_n$.

\smallskip
{\sc Remark.} As can be seen Theorems A and B are very similar.
However, in one sense Theorem B was easier to prove by induction
on $n$ because we could use the natural inclusions
(\ref{inclusionP}), whereas the natural inclusions of the real
projective plane were not transversal to the $G_n$-action and
had to be replaced by a second set of embeddings. On the other
hand, Theorem B has an exceptional fibre $P_n(\infty)$ for $n \geq 1$,
whereas in
Theorem A, $f_n$ is a fibration outside the branch locus. Notice
that, in each case, our branch locus $P_n$ is embedded in a
manifold of twice its dimension, namely
$$
\begin{array}{ll}
P_n \subset P_{n+1} & ({\rm Theorem\ A})\ ,\\
P_n \subset P_n(C) & ({\rm Theorem\ B})\ .
\end{array}
$$
These coincide only for $n=0$. For $n \geq 1$ even the normal
bundles are different.

\medskip
The four spheres $S^{d(n)}$ together with their cohomogeneity one
actions by $G_n$ appear in several contexts in differential
geometry. We discuss briefly two of them.

\smallskip
a) {\it Isoparametric hypersurfaces}. A real-valued function $f$
on a Riemannian manifold is isoparametric if the first
and second differential parameter of $f$ (i.e.\ ${\parallel}
\grad\,f {\parallel}^2$ and $\Delta f$) are constant along the
level sets of $f$. The interest in such functions originated from
geometrical optics. Any regular level set of an isoparametric
function is called an isoparametric hypersurface. E.\ Cartan
proved that a hypersurface in a space of constant curvature is
isoparametric if and only if it has constant principal curvatures.
The number of distinct principal curvatures of an isoparametric
hypersurface in a sphere is $1,2,3,4$ or $6$. It is easy to show
that the isoparametric hypersurfaces in spheres with $1$ or $2$
distinct principal curvatures are the distance spheres or
Clifford tori, respectively. In \cite{Car2} E.\ Cartan proved that
the isoparametric hypersurfaces with $3$ distinct principal
curvatures exist only in $S^{d(n)}$, $n=0,1,2,3$, and moreover
they are the regular orbits of the cohomogeneity one action of
$G_n$ on $S^{d(n)}$. For a survey about this topic see \cite{Tho}.

\smallskip
b) {\it Positive curvature}. A classical problem is to classify
all simply connected closed smooth manifolds which admit a
Riemannian metric with positive sectional curvature.
The standard examples of such
manifolds are the spheres $S^m$ and the projective spaces $CP^m$,
$HP^m$ and $OP^2$ ($m \geq 2$). Wallach \cite{Wal} proved that the
only simply connected closed smooth manifolds admitting a
{\it homogeneous}
Riemannian metric with positive sectional curvature are, apart
from the even-dimensional spheres and the above projective spaces,
precisely the regular orbits of the cohomogeneity one actions of
$G_n$ on $S^{d(n)}$, $n = 1,2,3$. Explicitly these are the flag
manifolds $SU(3)/U(1)^2$, $Sp(3)/Sp(1)^3$ and $F_4/Spin(8)$ of all
full flags in $CP^2$, $HP^2$ and $OP^2$. However, we should point
out that the positive curvature metrics on these flag manifolds
are not the induced metrics from the spheres with their standard
metrics.

\section{The Projective Version}

So far we have concentrated on constructing maps compatible with
the compact isometry groups $G_n$ of the projective planes $P_n$.
Now we shall extend these results appropriately to the non-compact
groups $\cG_n$ of projectivities of $P_n$.

In Section 6 we saw that the embeddings $P_n \subset H_n(1)$,
which are compatible with the action of $G_n$, extend naturally to
embeddings
\begin{equation*}
P_n \subset \cP_n
\end{equation*}
compatible with the action of $\cG_n$. Here $\cP_n$ is the real
projective space associated with $H_n$, and so is the projective
completion of the affine space $H_n(1)$. The vector space $H_n$ is
an irreducible representation of $\cG_n$ which splits off a
one-dimensional trivial factor (corresponding to
the trace) on restriction to $G_n$, so
that $\cG_n$ acts on $\cP_n$ with $G_n$ preserving the hyperplane
at infinity and the central point given by the scalar multiples of $I/3$.

The projection $\pi_n : H_n \to H_{n-1}$ induces a corresponding
projection
\begin{equation}
\pi_n : \cP_n \setminus \cP(H_{n-1}^{\perp}) \to \cP_{n-1}
\label{projection1}
\end{equation}
which is defined in the complement of the ``axis'' of projection
arising from $H_{n-1}^{\perp}$. The map (\ref{projection1}) is
compatible with the action of $\cG_{n-1}$. Note that $P_n$
is contained in $H_n(1)$ and so does not intersect the axis of the
projection (\ref{projection1}). Hence we get a well-defined map
\begin{equation}
\pi_n : P_n \to \cP_{n-1}
\label{projection2}
\end{equation}
compatible with the action of $\cG_{n-1}$. We propose to examine
(\ref{projection2}) with respect to the orbits of $\cG_{n-1}$.

We begin by looking at the action of $\cG_n$ on $\cP_n$ (and we
shall then replace $n$ by $n-1$). We recall that we have a cubic
polynomial $\det$ in $H_n$ whose vanishing defines a hypersurface
$Z_n \subset \cP_n$. In the appendix which follows, the
complexification $Z_n(C) \subset \cP_n(C)$ is discussed in detail.
The group $\cG_n$ leaves $Z_n$ invariant: this is clear for
$n=0,1$, requires a little verification for $n=2$,
and is a classical result of E.~Cartan \cite{Car} when
$n=3$ and $\cG_3 = E_6^{-26}$ is a real form of $E_6(C)$. Moreover $Z_n$
contains $P_n$ as a $\cG_n$-orbit and $P_n$ is the singular locus
of $Z_n$.

The group $G_n$ acts on the affine part $H_n(1)$ of $\cP_n$ with
$I/3$ as fixed point and its orbits are parametrized by three real
eigenvalues
\begin{equation*}
\lambda_1 \leq \lambda_2 \leq \lambda_3\ \ ,\ \ \lambda_1 +
\lambda_2 + \lambda_3 = 1\ .
\end{equation*}
$Z_n$ is given by $\lambda_1\lambda_2\lambda_3 = 0$, while $P_n$
is given by $\lambda_1 = \lambda_2 = 0$ and $\lambda_3 = 1$.
We can indicate the affine part of $Z_n$ schematically by the
following picture (which is actually what a $2$-dimensional (affine) slice
would look like: a real cubic curve):
\par
\begin{figure}[h]
\epsfbox{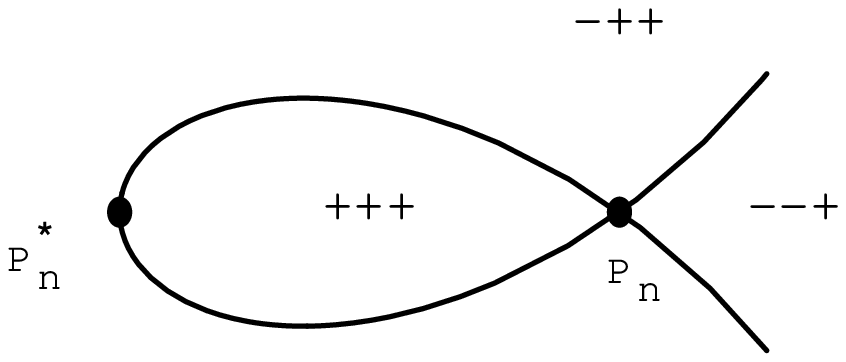}
\end{figure}
\noindent
where the complement of $Z_n$ in $H_n(1)$ is divided into
regions, depending on the signs of the three eigenvalues as
indicated. We
shall focus attention on the bounded region of positive-definite
matrices, which we denote by $\Lambda_n^+$,
and its boundary which will be denoted by $\Sigma_n$.
Notice that $\Sigma_n$ is a semi-algebraic set (given by
polynomial equations and inequalities), since it is only part of
the real algebraic variety $Z_n$.

Note that $\Sigma_n$ contains not only $P_n$, as a singular locus,
but also another copy $P_n^*$ given by
\begin{equation*}
\lambda_1 = 0\ ,\ \lambda_2 = \lambda_3 = \frac{1}{2}\ .
\end{equation*}
Unlike $P_n$ the copy $P_n^*$ consists of smooth points of $Z_n$.

Our first key lemma is

\medskip
{\sc Lemma 3.} {\it The set $\Lambda_n^+$ is convex and
its boundary $\Sigma_n$ is homeomorphic
to a sphere of dimension $d(n) = 3\cdot 2^n + 1$.}

\medskip
The convexity follows from the fact that $Z_n$ is a {\it cubic}
hypersurface so that any line meets it in at most three points.
Thus a chord of $\Sigma_n$ cannot exit and then reenter $\Lambda_n^+$.
Projection from any interior point, say
$I/3$, then gives the required homeomorphism with the sphere.

\medskip
{\sc Remark.} For the classical cases $n=0,1,2$ this lemma is
directly evident from the properties of eigenvalues of Hermitian
matrices: in fact, as we see, it also holds for the octonionic
case $n=3$. The classical cases were used by Arnold \cite{Arn} to
establish Theorem A for $n=1,2$. He also had results for larger
matrices. By contrast we stick to $3 \times 3$ matrices but handle
also the Cayley case.

\medskip
Our next lemma describes the $\cG_n$-orbit structure of $\cP_n$.

\medskip
{\sc Lemma 4.} {\it The $\cG_n$-orbits on $\cP_n$ are as follows:
\begin{itemize}
\item[(1)] Two open orbits, given by matrices in $\Lambda_n^+$ and by matrices
with $\lambda_1 < 0 < \lambda_2 \leq \lambda_3$;
\item[(2)] Two orbits of codimension one, namely $\Sigma_n
\setminus P_n$ (which is given by $\lambda_1 = 0 < \lambda_2 \leq \lambda_3$)
and $Z_n \setminus \Sigma_n$ (which is
given by $\lambda_1 < 0 = \lambda_2 < \lambda_3$);
\item[(3)] One orbit of codimension $2^n+2$, namely $P_n$,
given by $\lambda_1 = \lambda_2 = 0$ and $\lambda_3 = 1$.
\end{itemize}
}

\medskip
{\it Proof.} This follows from the fact that the $\cG_n$-orbits
on $H_n$ are characterized by the
rank of the matrix and the sign of the eigenvalues
(in the projective space $\cP_n$ a matrix $X$ is also equivalent to
$-X$).

\medskip
{Finally} we shall need to know the $\cG_{n-1}$-orbits on $P_n$:

\medskip
{\sc Lemma 5.} {\it The action of $\cG_{n-1}$ on $P_n$ has two
orbits, namely $P_{n-1}$ and its complement.}

\medskip
{\it Proof.} First we prove this for $n=1$. We have to show that
$SL(3,R)$ acts transitively on $CP^2 \setminus RP^2$. Let $\xi$ be
in this open set. Then $\xi \neq \bar{\xi}$, so we have a unique
$CP^1$ joining $\xi$ and $\bar{\xi}$. This $CP^1$ is
the complexification of a suitable line $RP^1$ in $RP^2$.
The subgroup of $SL(3,R)$ that
leaves this $RP^1$ invariant induces on $CP^1$
an action of $GL(2,R)$. Now $SL(2,R)$
acting on $CP^1 = S^2$ acts transitively on each hemisphere (the
model of the hyperbolic plane), and an element of determinant $-1$
switches the two hemispheres. Together with the fact that
$SL(3,R)$ is transitive on lines in $RP^2$ this establishes the
lemma for $n=1$. To prove it for $n=2,3$ we use the second
embeddings introduced in Section 5,
$$
CP^2(j) \subset HP^2\ \ {\rm and}\ \ g(CP^2) \subset OP^2\ ,
$$
which according to (\ref{intersectH}) and (\ref{intersectO})
intersect the original $CP^2$ in $RP^2$. We recall that the orbits
of the compact groups $Sp(3)$ and $F_4$ cut out on $CP^2(j)$ and
$g(CP^2)$ the orbits of $SO(3)$. Since we have shown that
$SL(3,R)$ acts transitively on $CP^2 \setminus RP^2$, and since
$$
SL(3,R) = \cG_0 \subset \cG_1 \subset \cG_2\ ,
$$
Lemma 5 follows for $n=2,3$.

\medskip
By definition the map $\pi_n : P_n \to
\cP_{n-1}$ of (\ref{projection2}) is the identity on
$P_{n-1}$, and
from Lemma 4 and Lemma 5 we see that $\pi_n$
must map $P_n$ into $\Sigma_{n-1}$, since this is
the only compact union of two $\cG_{n-1}$-orbits which lies in the
affine part $H_{n-1}(1)$ of $\cP_{n-1}$. It is easy to check that
we cannot have $\pi_n(P_n) = P_{n-1}$: it is enough to check this
for $n=1$, when it cannot happen for dimension reasons
($\pi_1(P_1)$ being of dimension $4$). Thus we deduce the
following projective refinement of Theorem A:

\medskip
{\sc Theorem A$^\prime$} {\it The map $\pi_n : P_n \to \cP_{n-1}$ has
image $\Sigma_{n-1}$ and $P_n \setminus P_{n-1} \to \Sigma_{n-1}
\setminus P_{n-1}$ is a homogeneous fibration for the group
$\cG_{n-1}$.}

\medskip
{\sc Remarks.} 1) Note that $\Sigma_{n-1}$ does not have a natural
smooth structure compatible with the action of $\cG_{n-1}$, since
$P_{n-1}$ is a singular locus. However, if we restrict to the
compact subgroup $G_{n-1} \subset \cG_{n-1}$, which fixes $I/3$,
then projection from $I/3$ maps $\Sigma_{n-1}$ to the sphere
$S^{d(n-1)}$ (as shown by Lemma 3) and now the map
$$
P_n \to S^{d(n-1)}
$$
identifies the smooth structure of $S^{d(n-1)}$ with the quotient
smooth structure of $P_n$ as explained in Section 2. Thus Theorem
A follows from Theorem A$^\prime$ and Lemma 3. The smoothing of
$\Sigma_{n-1}$ by the radial projection effectively ``rounds off
the corners'' as exemplified by the projection of a square from
its centre onto a circle.

2) By considering the projective lines of $P_n$ we get the
following more precise picture of the geometry of $Z_n$ in
relation to $P_n$. Consider any line in $P_n$. This is a sphere of
dimension $2^n$ embedded in the standard way as a real quadric in
$\cP_n$ (lying in a linear subspace of dimension $2^n + 1$). Take
its interior, an open ball, and its closure. This lies inside
$\Sigma_n$, and $\Sigma_n$ is filled up by the union of the closed
balls. If we fix a standard metric on $P_n$, and hence a compact
subgroup $G_n \subset \cG_n$, then the sphere becomes a round
sphere in $H_n(1)$ and its centre lies on $P_n^*$, which is the
locus of such centres under the action of $G_n$. This shows that
$P_n^*$ is the closest part of $\Sigma_n$ to the centre $I/3$. For
a fixed line in $P_n$ the geometry of its interior is just
hyperbolic geometry and the subgroup of $\cG_n$ preserving the
line is the conformal group of the sphere or the isometry group of
its interior. The corresponding complex picture is described in
the appendix.

3) Since $\pi_n$ is a linear projection it takes the interiors of
real quadrics to the interiors of their projections. Remark 2 then
shows that $\pi_n$ maps $\Sigma_n \setminus P_n$ into the interior
of $\Sigma_{n-1}$. Finally, since this latter is an open orbit of
$\cG_{n-1}$, it follows that we get the whole of the interior.

\medskip
We shall now consider to what extent there is a Theorem B$^\prime$
analogous to Theorem A$^\prime$. It is easy to see that there can
be no version compatible with $\cG_{n-1}$, since this moves the
projective plane $P_{n-1}^*$ while, for Theorem B (for $n
> 0$), $P_{n-1}^*$ is a distinguished subspace of the sphere,
being the image of the exceptional fibre. However the first part
of Theorem A$^\prime$ extends to give

\medskip
{\sc Theorem B$^\prime$.} {\it The image of the map $\sigma_n :
P_n(C) \to H_n(0)$ is $\Sigma_n$.}

\medskip
{\it Proof.} For $n=0$ we have $\sigma_0 = \pi_1$ (of Theorem A$^\prime$)
and hence Theorem B$^\prime$ for $n=0$ follows from Theorem A$^\prime$ for
$n=1$. For $n \geq 1$ we have a commutative diagram
$$
\CD
P_{n-1}(C) @>\sigma_{n-1}>> H_{n-1}(0)   \\
@V VV @VV  V   \\
P_n(C)  @>>\sigma_n> H_n(0)
\endCD \ \ \ \ .
$$
Moreover the codimension one orbits of $G_n$ in $P_n(C)$ cut out
the codimension one orbits of $G_{n-1}$ on $P_{n-1}(C)$. Hence the
image of $\sigma_n$ is the union of the $G_n$-orbits of the image
of $\sigma_{n-1}$. But, from the characterization of the
$G_n$-orbits in $H_n$ by their eigenvalues, it then follows that
$$
G_n(\Sigma_{n-1}) = \Sigma_n\ .
$$
Thus Theorem B$^\prime$ follows by induction on $n$.

\section{Appendix}

In Section 6 we studied the differential geometry of the
complexified projective varieties $P_n(C)$. In this appendix we will
review (without complete proofs) some of the algebraic geometry
which is of independent interest and provides further background.
All dimensions in this section are complex dimensions. For further
details see \cite{LaMa}.

The algebraic variety
$$
P_n(C) \subset \cP_n(C)\ \ \ (n = 0,1,2,3)
$$
of dimension $2^{n+1}$ has appeared as an orbit of the
complex Lie group $\cG_n(C)$ in the complex
projective space $\cP_n(C)$ of dimension $3\cdot 2^n + 2$.
The group $\cG_n(C)$ is the complexification of the
compact Lie group $\hat{G}_n$ given in table (\ref{table}). Explicitly
we have the sequence of groups $\cG_n(C)$ ($n=0,1,2,3$)
\begin{equation*}
SL(3,C)\ ,\ SL(3,C) \times SL(3,C)\ ,\ SL(6,C)\ ,\ E_6(C)\ .
%\label{complexgroups}
\end{equation*}
These have irreducible representations on the vector spaces
$H_n(C) = H_n \otimes C$, which are the complexifications of the
real vector spaces $H_n$ of all $3\times 3$
Hermitian matrices over the division algebra $A_n$.
These representations have dimension $3(2^n+1)$, explicitly (for $n=0,1,2,3$)
$$
6 \ ,\ 9\ ,\ 15\ ,\ 27\ ,
$$
and they projectivize to give the spaces $\cP_n(C)$.

We can also consider the dual representations on $H_n(C)^*$ giving
dual projective spaces $\cP_n(C)^*$, and in these there is a unique
compact orbit $P_n(C)^*$ of $\cG_n(C)$.

On $H_n(C)$ there is a unique (up to scalars) cubic polynomial
invariant under $\cG_n(C)$ which we have denoted by $\det$, for
reasons explained in Section 3. For $n = 0,1$ it is just the usual
determinant. For $n=2$ it is the $SL(6,C)$-invariant cubic on
$\Lambda^2(C^6)$ given by the exterior cube into $\Lambda^6(C^6)
\cong C$.
For $n=3$ with $\cG_3(C) = E_6(C) \subset SL(27,C)$ this invariant
cubic was discovered by E.\ Cartan \cite{Car}. In all cases it defines a cubic
hypersurface in $\cP_n(C)$ which we denote by $Z_n(C)$, and which
contains $P_n(C)$.

The action of $\cG_n(C)$ on $\cP_n(C)$ has just three orbits,
namely $P_n(C)$, $Z_n(C) \setminus P_n(C)$ and $\cP_n(C) \setminus
Z_n(C)$. When $n = 0$ these just correspond to symmetric matrices
of ranks $1,2,3$, and we could use the same terminology in the
general case. In fact $Z_n(C) \setminus P_n(C)$ is the
$\cG_n(C)$-orbit of the diagonal matrix $\Diag(1,1,0)$ of rank $2$, and
as we have already observed $P_n(C)$ is the $\cG_n(C)$-orbit of
the diagonal matrix $\Diag(1,0,0)$. Points of rank $3$ constitute the
$\cG_n(C)$-orbit of the unit matrix. The complex Lie group $\cG_n(C)$ is just
the identity component of the group of holomorphic transformations
of $P_n(C)$.

The natural embedding of $2 \times 2$ matrices into $3 \times 3$
matrices by adding zeroes in the third row and column
gives a linear subspace
\begin{equation*}
\cL_n(C) \subset \cP_n(C)\ ,\ \dim \cL_n(C) = 2^n+1\ ,
%\label{linsub}
\end{equation*}
and its orbit under $\cG_n(C)$ fills out the whole of $Z_n(C)$.
The intersection
\begin{equation*}
\cL_n(C) \cap P_n(C) = L_n(C)
%\label{intersect}
\end{equation*}
is a complex quadric of dimension $2^n$. In fact, $L_n(C)$ inside
$\cL_n(C)$ plays the same role (for $2 \times 2$ matrices) that
$Z_n(C)$ does inside $\cP_n(C)$: it is the hypersurface given by
the invariant {\it quadratic} $\det_2$ for $2 \times 2$ matrices.

The family of all transforms of $\cL_n(C)$ under $\cG_n(C)$
therefore cuts out on $P_n(C)$ a corresponding family of quadrics
$L_n(C)$ with
$$
\dim L_n(C) = \frac{1}{2}\dim P_n(C) \ .
$$
The dual $\cL_n(C)^*$ of $\cL_n(C)$ is a linear subspace of
$\cP_n(C)^*$ of dimension
$$
(3 \cdot 2^n + 2) - (2^n + 1) - 1 = 2^{n+1} = \dim P_n(C)^*\ .
$$
In fact, $\cL_n(C)^*$ is the tangent space to $P_n(C)^*$ at a
point $\ell_n$. The correspondence
$$
L_n(C) \longleftrightarrow \ell_n
$$
represents $P_n(C)^*$ as the parameter family of the quadrics
$L_n(C)$ on $P_n(C)$, and the situation is symmetrical (or dual):
points of $P_n(C)$ parametrize quadrics (of half the dimension)
on $P_n(C)^*$.

Note that the quadrics $L_n(C)$ are all non-singular (since all
points are of ``rank 1''), whereas $Z_n(C)$ has $P_n(C)$ as a
{\it singular locus}: a generic point of $Z_n(C)$ has rank 2,
whereas points of $P_n(C)$ are of rank 1. Thus $Z_n(C)$ determines
$P_n(C)$ (as its singular locus), and conversely $P_n(C)$
determines $Z_n(C)$, as the space generated by the linear
spaces $\cL_n(C)$ spanned by the quadrics $L_n(C)$.

It may be helpful at this stage if we looked in detail at the
special case $n=0$, so that we are dealing with the classical
embedding of $CP^2$ as the Veronese surface $V$ in $CP^5$. The
lines of $CP^2$ become conics on $V$ and these lie in planes. The
lines of $CP^2$ are parametrized by the dual $CP^2$ which can be
identified with the dual Veronese surface $V^*$ in $(CP^5)^*$.
Thus the planes spanned by the conics in $V$ form a $2$-parameter
family and they fill out a (cubic) hypersurface $Z$. Since every
pair of distinct points on $CP^2$ lies on a unique line, every
pair of distinct points of $V$ lies on a unique conic. In
particular it follows that the chordal variety of $V$ (i.e.\ the
closed subspace generated by all chords) is also the space
generated by all the planes spanned by the conics and hence is the
hypersurface $Z$. This is a very unusual situation for a surface
in $CP^5$. On dimension grounds one could expect the chordal
variety to be the whole ambient space. Equivalently, when the
chordal variety is only a hypersurface, the projection from a
generic point gives an {\it embedding}\footnote{Over the reals
this gives an embedding $RP^2 \subset RP^4$. Since this is of
degree 2 it lifts to the double cover $S^4$ giving the embedding
$RP^2 \subset S^4$ of (\ref{branchloci}).}
(without singularities) in
$CP^4$. In fact, it is a classical result of Severi \cite{Sev} that
the Veronese surface is the only surface (not contained in
a hyperplane) in $CP^5$ with this property.

Zak \cite{Zak} (see also \cite{LaVa})
has investigated this ``Severi property'' for
higher dimensions when $V_d \subset CP^N$. The critical case is
when $d = \frac{2}{3}(N-2)$. Zak proved the remarkable result
that, firstly
$$
d = 2^{n+1}\ ,\ n=0,1,2,3\ ,
$$
and secondly that the only such varieties in these dimensions
are the complexified projective planes $P_n(C)$ in their
standard projective embeddings in $\cP_n(C)$. For this reason
these varieties have been named Severi varieties \cite{LaMa}.

To see how this fits into the picture we have described we have
to note that in all cases {\it the chordal variety of $P_n(C)$ is
the cubic hypersurface $Z_n(C)$}. The proof is very similar to the
case $n=0$, but with a caveat. Given any two distinct points $x$
and $y$ in $P_n(C)$ there are just two possibilities, either
\begin{itemize}
\item[(i)]  the projective line in $\cP_n(C)$
containing $x$ and $y$ lies entirely on
$P_n(C)$, or
\item[(ii)] there is a unique quadric $L_n(C)$ on $P_n(C)$ containing both $x$
and $y$.
\end{itemize}
For $n = 0$ case (i) never happens: the classical Veronese variety
contains no lines. Clearly (ii) is the generic situation and as
before this implies that the {\it hypersurface $Z_n(C)$} (generated by the
planes $\cL_n(C)$ spanned by the $L_n(C)$) {\it is precisely the
chordal variety of $P_n(C)$.} This shows that $P_n(C)$ does indeed
have the Severi property in $\cP_n(C)$.
The power of Zak's Theorem is that these are the
only ones.

There is a striking resemblance between Zak's theorem in complex
algebraic geometry and the classical results about division
algebras and projective planes. It would be interesting to see if
a purely topological proof of Zak's theorem could be found. We
recall that the use of Steenrod squares enables one to prove that
a projective plane must have dimension a power of 2 (analogous to
the first part of Zak's Theorem), while K-theory is needed for the
final part \cite{AdAt}. One is therefore tempted to expect a K-theory proof of
Zak's theorem, particularly in view of the role that K-theory, as
developed by Grothendieck, plays in algebraic geometry.

It is at this point that we should perhaps pass from the purely
complex approach to the varieties $P_n(C)$ and introduce real
structures. Recall that $P_n(C)$ has a complex conjugation,
preserved only by the real subgroup $\cG_n$ of $\cG_n(C)$, and
that the set of real points just recover the original projective
planes $P_n$. The complex quadrics $L_n(C)$ on $P_n(C)$,
parametrized by $P_n(C)^*$, include those preserved by conjugation
and parametrized by $P_n^*$. These are just the complexifications
of the projective lines in $P_n$ (i.e.\ spheres of dimension
$2^n$). The incidence properties of these projective lines in
$P_n$ imply generically the corresponding properties of the
complex quadrics $L_n(C)$ in $P_n(C)$, ensuring (ii) above.
However, for $n \geq 1$, the exceptional case (i) does occur. For
example, when $n =1 $, we have $P_1(C) = CP^2 \times CP^2$, and two
distinct point pairs whose first components agree give case (i).

The manifolds $P_n(C)$, with their family of submanifolds $L_n(C)$
of middle dimension, are examples of generalized projective planes
in the sense of Atsuyama \cite{Ats}, who studied these from the
point of view of differential geometry.

Recall that we have the following subgroups of $\cG_n(C)$:
$$
\CD
G_n @> >> \cG_n   \\
@V VV @VV  V   \\
\hat{G}_n  @>> > \cG_n(C)
\endCD \ \ \ \ ,
$$
where the groups in the first column are the maximal compact
subgroups of those in the second column and preserve the metrics on
$P_n$ and $P_n(C)$ respectively. We can also introduce the groups
$G_n(C)$, the complexification of $G_n$. Since the representation
$H_n(C)$ splits off a trivial factor of dimension one (given by the
trace) when restricted to $G_n(C)$, it follows that there is a
hyperplane section $P_n(\infty)$ of $P_n$ invariant under
$G_n(C)$. In fact, $G_n(C)$ acts on $P_n(C)$ with just two orbits,
$P_n(\infty)$ and its complement. We propose to examine the
geometry of $P_n(\infty)$ which is a homogeneous space of $G_n(C)$
and so also of $G_n$. In particular we will see that we can
reconstruct $P_n(C)$ canonically from $P_n(\infty)$.

It is instructive to consider first the simple case $n=0$. Then
$P_0(\infty)$ is a rational normal quartic curve, the image under
the Veronese embedding of a conic in $CP^2$. This is the conic
$z_1^2 + z_2^2 + z_3^2 = 0$ invariant under $G_0(C) = SO(3,C)$.
Consider now the conics on $V = P_0(C)$ (the Veronese surface)
which are complexifications of projective lines in $P_0 = RP^2$.
In the ``abstract'' $CP^2$ which maps to $V$ these just correspond
to projective lines with real equations. Each such line meets the
conic $z_1^2 + z_2^2 + z_3^2 = 0$ in a pair of conjugate points.
The quotient of the conic by this involution is naturally
identified with the {\it dual} $RP^2$. This is part of the content
of Theorem B for $n=0$. If we now consider all conics on $V$, they
come from all projective lines on $CP^2$, and they meet the conic
$z_1^2 + z_2^2 + z_3^2 = 0$ in any pair of points. Thus if we
start with our rational normal curve $P_0(\infty)$ on $P_0(C)$ and
consider the variety of all (unordered) pairs of points on
$P_0(\infty)$ (i.e.\ its symmetric square) we get the dual $CP^2$.
From this we can (by duality) recover the original $CP^2$ and also
its Veronese embedding. Thus $P_0(\infty)$ determines the whole
picture.

We want to show that this is typical of the general case (i.e.\
for all $n$).

We consider the real family of complex quadrics $L_n(C)$
parametrized by points of $P_n^* \subset P_n(C)^*$. These are just
the complexifications of the projective lines of the projective
plane $P_n$. Since any two distinct points of $P_n$ are joined by
a unique line they are never in the special position of (i).
Dually this means that no two quadrics of the real family of
$L_n(C)$ meet in more than one point. Since they already have one
common point on $P_n$ they meet nowhere else. In particular the
family of quadrics of one lower dimension cut out on
$P_n(\infty)$,
$$
L_n(\infty) = L_n(C) \cap P_n(\infty)\ ,
$$
are all disjoint. A dimension count shows that they must fill out
the whole of $P_n(\infty)$, and since they are by construction
parametrized by the dual space $P_n^*$ we get a fibration
\begin{equation}
P_n(\infty) \to P_n^*
\label{fibrationPninfty}
\end{equation}
with fibre $L_n(\infty)$. This gives a more explicit description
of the behaviour at infinity of the map of Theorem B.
Note that, for $n = 0$, $L_0(\infty)$ has dimension zero and is
a point-pair as we have already seen.

Motivated by the case $n = 0$ we now consider the full complex
family of quadrics $L_n(C)$. These intersect $P_n(\infty)$ in the
full complex family of quadrics (of one lower dimension) whose
real members are the fibres of (\ref{fibrationPninfty}). This full
family is therefore parametrized by the complexification
$P_n(C)^*$ of $P_n^*$. Since we have cut down the symmetry from
$\cG_n(C)$ to $G_n(C)$ there is a distinguished subset
$P_n(\infty)^*$. These quadrics in $P_n(\infty)$ are {\it
singular}, they arise from quadrics $L_n(C)$ which touch
$P_n(\infty)$: note that these are never the real members (the
fibres of (\ref{fibrationPninfty})) since $P_n^*$ and
$P_n(\infty)^*$ are disjoint. All of this checks with what we saw
for $n = 0$.
Thus $P_n(\infty)$ contains a family of complex quadrics
$L_n(\infty)$, parametrized by $P_n(C)^*$. This enables us to
recover $P_n(C)^*$ and hence $P_n(C)$ from $P_n(\infty)$.

The fibration (\ref{fibrationPninfty}) can be viewed as a {\it
twistor fibration}. For example when $n = 1$, $P_1(C) = CP^2
\times (CP^2)^*$ and $P_2(\infty)$ is the incidence locus. It is
therefore the flag manifold of $SU(3)$ and
(\ref{fibrationPninfty}) is the twistor fibration for $CP^2$
regarded as a $4$-manifold with self-dual metric \cite{AtHiSi}.
For $n=2,3$ the fibration (\ref{fibrationPninfty}) is a partial
twistor fibration in the sense of Bryant \cite{Bry}, and it is
given as an interesting example of a general theory. Explicitly,
the two fibrations are
$$
\CD P_2(\infty) = \frac{\displaystyle Sp(3)}{\displaystyle Sp(1)U(2)}
@>>>\frac{\displaystyle Sp(3)}{\displaystyle Sp(1)Sp(2)} = HP^2 \cong
P_2^*
\endCD
$$
and
$$
\CD P_3(\infty) = \frac{\displaystyle F_4}{\displaystyle Spin(7)U(1)}
@>>>
\frac{\displaystyle F_4}{\displaystyle Spin(9)} = OP^2 \cong P_3^* \
,
\endCD
$$

\smallskip\noindent
where the fibres are $L_2(\infty) = Sp(2)/U(2)$
and $L_3(\infty) = Spin(9)/Spin(7)U(1)$ respectively.
Note that $Sp(2)/U(2)$ is isomorphic to the $3$-dimensional quadric
$SO(5)/SO(3)SO(2)$ and that $Spin(9)/Spin(7)U(1)$ is isomorphic to
the $7$-dimensional quadric $SO(9)/SO(7)SO(2)$.

\end{document}